\documentclass[11pt]{article}
  \usepackage[utf8]{inputenc}
\usepackage{amsmath,amssymb,epsf}
\usepackage{amsmath}
\usepackage{graphics}
\advance \textwidth 3cm
\advance \textheight 3cm
\usepackage{amsfonts}
\usepackage{latexsym}
\newtheorem{theoreme}{Theorem}
\newtheorem{lemme}{Lemma}

\newtheorem{definition}{Definition}

\newtheorem{prop}{Property}
\newtheorem{remarque}{Remark}
\newenvironment{preuve}[1]{\par\noindent\underline{Proof #1} :\quad}%
    {\unskip\nobreak\hfil\penalty50\hskip2em\null\nobreak\hfil%=
$\Box$\parfillskip0pt\par\medskip}%versions F et A

\newcommand{\spec}{\mathrm{Spec}}

\newcommand{\vect}{\mathrm{vect}}

\title{ Asymptotic behavior of the eigenvalues of Toeplitz matrices with even symbol}

\author{ Philippe Rambour\thanks{Universit\'{e} de Paris Saclay,
      B\^atiment 307; F-91405
Orsay Cedex;
tel : 01 69 15 57 28 ; 
      \mbox{e-mail : philippe.rambour@universite-paris-saclay.fr}}
       }
\date{}
\begin{document}
\maketitle
            \renewcommand{\abstractname}{Abstract}
          \begin{abstract}
In this paper we consider  an interval $[\theta_{1}, \theta_{2}] \subset ]0, \pi[$ 
and $f$ a periodic and even function in $C^4 \left([0, 2 \pi]\right)$  such that 
  $f(\theta) \in [f(\theta_{1}), f(\theta_{2})] \iff \theta \in [\theta_{1}, \theta_{2}]$ and
    $f\in \mathcal A^+\left([\theta_{1}, \theta_{2}]\right)$
   (resp. $f\in \mathcal A^-\left([\theta_{1}, \theta_{2}]\right)$, 
  that means  $f'(\theta) >0$ (resp. $f'(\theta) <0$) for all 
   $\theta \in [\theta_{1}, \theta_{2}]$. 
Then we obtain a higher order asymptotic formula for all  the eigenvalues of the Toeplitz matrix $T_N(f)$ as $N \to +\infty$ which 
belong to $[f(\theta_{1}), f(\theta_{2})]$ (resp. $[f(\theta_{2}), f(\theta_{1})])$.     
 \end{abstract}

   \vspace{1cm}


\vspace{1cm}


\section{Introduction and statement of the main results}
If $\mathbb T = R / (2\pi \mathbb Z)$ and  $h\in L^{1} (\mathbb T)$, we denote by $T_{N}(h)$ the Toeplitz matrix of order $N$ with symbol $h$. It is  the $(N+1)\times (N+1)$ matrix such that, for $N\ge k,l\ge 0$,
$\left( T_{N}(h)\right)_{k+1,l+1}= \hat h (k-l)$ where $\hat h(u)$ is the Fourier coefficient of order $u$ of $h$ (\cite{GS,Bo.3}). For a real valued function $h$ the matrix  $T_{N}(h)$ is a Hermitian Toeplitz matrix. 
We here consider a symmetric Toeplitz matrix, which is equivalent to assuming that the symbol $h$ is an even function and  we denote by 
$ \lambda_{N}^{(1)} \le \lambda_{N}^{(2)}\le \cdots \le  \lambda_{N}^{(N+1)}$ the eigenvalues of $T_{N}(h)$. 
This paper adresses the asymptotic behavior of the eigenvalues of $T_{N}(h)$ as $N$ goes to infinity. This is a topic which has attracted mathematicians and physicists for a long time.  Toeplitz matrices 
and their relatives emerge in particular in statistic \cite{Dahlhaus,B} and in statistical physics \cite{Bas2,Bas3}. It is known from a long time that Toeplitz matrices are useful for providing Green's kernels for studying the solutions of certain differential equations (\cite{SpSt,RS04}), and also for discretizing differential operators with finite differences \cite{STMN1}. Finally, these matrices are used in more recent fields, such as Ising models \cite{DeiftKra2} and iso-geometric analysis \cite{Garoni1}. The questions about the asymptotic behavior of their spectral characteristics, especially their determinants, eigenvalues, and eigenvectors, are always at the heart of the matter. We refer to the papers \cite{DeiftKra2} for an extensive list of references. According to  the first Szeg\H{o} limit theorem (see \cite{GS}) the eigenvalues of 
$T_{N}(h)$ are asymptotically distributed as the value of $h$; see \cite{GS} for $L^{\infty}$ symbols, 
\cite{Tyr2} for $L^1$ symbols, and \cite{TRENCH2,Tyr} for more general situations. 
In the Hermitian case extensive works has been done on the search for eigenvalues (or the extreme eigenvalues) of Toeplitz matrices \cite{Wid,GS,Serra,JMR01,Part,Part2,Part3} and more recently, for instance,  
\cite{BBGM2,BoBOGru,BabOGruMax,ESR,EST2,BoGru23}.In \cite{BBGM2} the authors give an asymptotic expansion  of order 2 for the eigenvalues of a Toeplitz matrice with smooth simple-loop symbol, the results of this article are the closest to ours, but the techniques are different.
 \cite{BoBOGru}  is a good reminder of the various results obtained on the eigenvalues of Toeplitz matrices with polynomial symbol. In 
 \cite {BabOGruMax} M. Barrera, A. B\H{o}ttcher, S. M. Grudsky and E. A. Maximenko show that the eigenvalues of the matrix $T_N (4(1-\cos\theta)^2)$ cannot have an asymptotic expansion to order 4. For banded Toeplitz matrices or block symmetric Toeplitz matrices, the reader is referred to \cite{ESR,EST2}. Here we consider real symbol. For complex symbols an alternative to our results is given in  \cite{DeiftITS}\\
The results of Theorem \ref{THEO1ter} are consistent with those of Theorem 2.3 of \cite{BBGM2}. But the method of proof is different and the statement concern functions which are outside the framework of \cite{BBGM2}. On the other hand  Theorem \ref{THEO1ter} indicates that the problem of the eigenvalues of Toeplitz matrices is a local problem, related to the variation of the function which is the symbol of the matrix. \\  %%%%
Here we denote by $\mathcal A^{+}([0,2\pi])$ (resp. $\mathcal A^{-}([0,2 \pi])$) the set of even differentiable periodic functions of period $2\pi$, 
such that $f'(\theta)>0$ (resp. $f'(\theta)<0$) for all $\theta$ in $]0,\pi[$. \\ 
More generally $[\theta_{1}, \theta_{2}] \subset ]0, \pi[$ we say that 
  $f\in \mathcal A^+\left([\theta_{1}, \theta_{2}]\right)$
   (resp. $f\in \mathcal A^-\left([\theta_{1}, \theta_{2}]\right)$ 
   if  $f$ is a differentiable, $2 \pi$ periodic, and even function such that $f'(\theta) >0$ (resp. $f'(\theta) <0$) for all 
   $\theta \in [\theta_{1}, \theta_{2}]$  
   and also $f(\theta) \in [f(\theta_{1}), f(\theta_{2})] \iff \theta \in [\theta_{1}, \theta_{2}]$ (resp.  $f(\theta) \in [f(\theta_{2}), f(\theta_{1})] \iff \theta \in [\theta_{1}, \theta_{2}]$).\\
For $\nu\ge0$ we denote by $W^\nu$ the weighted Wiener algebra 
of all functions $\psi : \mathbb T \mapsto \mathbb C$ which admits the representation 
$\psi  (t) = \displaystyle{\sum_{j\in \mathbb Z} \hat \psi (j) t^j}$ whose Fourier coefficients satisfy
$$ \Vert \psi \Vert_{\nu}=\sum_{j \in \mathbb Z} \vert \hat \psi (j) \vert \left( \vert j \vert +1\right)^\nu < \infty.$$
  Now for $f \in \mathcal A^+\left([\theta_{1}, \theta_{2}]\right)$
we define the functions 
   \[ H(\theta', \theta^{\prime\prime}) = 
  \frac{f(\theta') - f(\theta^{\prime\prime})}{(1-\cos\theta') - (1-\cos\theta^{\prime\prime})}  \quad \mathrm{for} \quad  (\theta' , \theta^{\prime \prime})\in [0, \pi]^2 \] 
  and 

  $$ \rho (\theta) = 
     \frac{1}{4\pi} P.V. \int_{0}^{2\pi}\frac{ \ln\left((H (t,\theta)\right)}{\tan \left( \frac{t-\theta}{2}\right)} dt 
    - \frac{1}{4\pi} P.V. \int_{0}^{2\pi}\frac{ \ln\left((H (t,\theta)\right)}{\tan \left( \frac{t+\theta}{2}\right)} dt \quad \mathrm{for} \quad \theta\in [\theta_{1},\theta_{2}].$$
    Lastly for an interval $[a,b] \subset [\theta_{1},\theta_{2}]$
 and an integer $N$ we denote by $k_{a,N}$ and $k_{b,N}$
the integers such that :
$k_{a,N}=\min\{k \vert \frac{k\pi}{N+2} \in [a,b]\}$,
$k_{b,N}=\max\{k \vert \frac{k\pi}{N+2} \in [a,b]\}$.  
We have also to define the two functions 
$$c_{1} (t) = f'(t) \rho(t)$$
and 
$$ c_{2} (t) =  f'(t) \rho'(t) \rho(t)+ \frac {1}{2} f^{\prime \prime} (t) \rho(t)^2;$$
Now we can state our main result and an easy consequence.
\begin{theoreme}\label{THEO1ter}
 Let $f$ in $C^4[0, 2\pi]$ be such that $f \in \mathcal A^+\left([\theta_{1}, \theta_{2}]\right)$ for an interval $[\theta_{1}, \theta_{2}] \subset]0, \pi[$.
Then for all interval $[a,b]\subset ]\theta_{1}, \theta_{2}[$ and for a  sufficiently large integer $N$ we have the two following statements.
\begin {enumerate}
\item
For all eigenvalue $\lambda$ of $T_{N}(f)$ in $[f(a),f(b)]$ we have a single integer $k$ in $[k_{\theta_{1},N},k_{\theta_{1},N}]$
such that $\lambda = {\tilde\lambda_{N}}^{(k)} + O \left( \frac{1}{(N+2)^3}\right)$, uniformly in $\lambda$,  with
 $$ {\tilde\lambda_{N}}^{(k)}= f \left(\frac{k \pi}{N+2}\right) + \frac{c_{1} \left(\frac{k \pi}{N+2}\right)}{(N+2)}+
 \frac{c_{2} \left(\frac{k \pi}{N+2}\right)}{(N+2)^2}.$$
 \item
     For all $k \in [k_{a,N},k_{b,N}]$  the matrix $T_{N}(f)$ has a single
eigenvalue $\lambda$ in $[f(\theta_{1}),f(\theta_{2})]$  such that $\lambda= {\tilde\lambda_{N}}^{(k)} + O \left( \frac{1}{(N+2)^3}\right)$ uniformly
 in $k$.
 \end{enumerate}
 \end{theoreme} 
\begin{remarque}
Similar results to Theorem \ref {THEO1ter} holds for the case where $f$ in $\mathcal A^-\left([\theta_{1}, \theta_{2}]\right)$ for an interval 
$[\theta_{1}, \theta_{2}]\subset ]0, \pi[$.
 \end{remarque}
\begin{remarque}
 In Theorem \ref {THEO1ter}  the eigenvalue $\lambda = \tilde \lambda_{N}^{(k)}+O(\frac{1}{(N+2)^2})$ is not necessarily $ \lambda_{N}^{(k)} $.
 \end{remarque}
 \begin{remarque}
If we consider the functions  $\psi_{\alpha} : \theta \mapsto (1-\cos \theta)^\alpha c(\theta)$ where $\alpha\ge 2$ and $c$ a  even positive  function such that 
 $c \in  C^4[0, 2\pi]$ and $ \alpha \sin \theta (1-\cos \theta)^{\alpha-1}  c(\theta) +  (1-\cos \theta)^\alpha c'(\theta)>0$ for all $\theta\in [0,\pi]$ 
 we can remark that Theorem \ref{THEO1ter} provides all the eigenvalues of the functions  of  $\psi_{\alpha}$ in $[\epsilon, \pi-\epsilon]$ 
 for all  $\epsilon>0$, that  is an extension of the main results of \cite{BoGru23}.
 \end{remarque}
   \begin{remarque}\label{REM3}
 Revisiting the proof of Theorem \ref{THEO1ter}, one can show that under the assumption 
 $f$ in  $C^3 ([0, 2 \pi]) $ we obtain an analogous version of this Theorem  with the formula
 $ \lambda = {\hat \lambda}_{N}^{(k)} +  O\left(\frac{\log N}{(N+2)^2} \right),$ with
  $${\hat \lambda}_{N}^{(k)} = f \left(\frac{k\pi}{N+2}\right)+ \frac{c_{1}(\frac{k\pi}{N+2})}{(N+2)}
$$ 
where the function $c_{1}$ is as in Theorem \ref{THEO1ter}.
\end{remarque}
   
  Our result can also be compared with that of Trench
  \cite{TRENCH93} where it is proved that for this class of symbols the eigenvalues are all distinct.\\
  To conclude we can remark that a tiny modification of the proof of Theorem \ref {THEO1ter} allows us to obtain the following Theorem
which is in fact Theorem 2-3 in \cite {BabOGruMax}.
\begin{theoreme}\label{THEO3}
 Let $f$ in $C^4[0, 2\pi]$ be such that $f \in \mathcal A^+\left([0,2 \pi]\right)$, $f^{\prime \prime} (0) > 0$ and $f^{\prime \prime} (\pi) < 0$, then for a  sufficiently large $N$ we have  
 $$ \lambda_{N}^{(k)}= f\left(\frac{k \pi}{N+2}\right) + \frac{c_{1} \left(\frac{k \pi}{N+2}\right)}{(N+2)}+
 \frac{c_{2} \left(\frac{k \pi}{N+2}\right)}{(N+2)^2}+O\left( \frac{1}{(N+2)^3}\right)$$
uniformly in $k=1, 2, \cdots, N$ and with $c_{1}$ and $c_{2}$ defined as previously.
\end{theoreme}
%%%%
  Lastly we have to recall the two following definitions 
 \begin{definition} 
 We denote by $\mathbb H^{+}$ is the set of all functions $\varphi$ in  
   $ L^2 (\mathbb T)$ whose Fourier coefficients satisfy   
  $ \hat \varphi (j) =0$  for all  $j <0$. $\mathbb H^{+}$ is called the Hardy space on the unit circle.
 \end{definition}
\begin{definition} 
 In the rest of this paper we denote by $\chi$  the function $\theta \mapsto e^{i\theta}$.
  \end{definition}
    \section{Proof of Theorem \ref{THEO1ter}}
   \subsection{Preliminaries}
   In this proof we have to use the following Theorem which 
   provides an inversion formula for a family of Toeplitz matrices.
\begin{theoreme} \label{THEO2} Let $P_{N+1}$ a trigonometric polynomial with degree $N+1$ 
and without zeros on the unit disc  $\bar D$.
Let $\omega = r \bar \chi_{0}$, $0<r<1$, $\vert \chi_{0} \vert =1$ and also 
$f_{r}=g_{1}g_{2,}$ with $g_{1}= \chi_{0} (1-\omega \chi) (P_{N+1})^{-1}$
and $g_{2}= (1-\omega \bar\chi) \left(\overline{P_{N+1}}\right)^{-1}$.
 Then for all polynomial $P$ in $\mathcal P_{N} = \vect \{1, \chi, \cdots , \chi^N\} $ we have 
 $$ T_{N}(f_{r})^{-1} (P) =
 \frac{1}{g_{1}} \pi_{+}\left( \frac{P}{g_{2}} \right)
 - \frac{1}{g_{1}} \pi_{+}\left( \Phi_{N}\sum_{s=0}^{+\infty} \left( H^{*}_{\Phi_{N}}
  H_{\Phi_{N}}\right)^{s} \pi_{+} \left( \tilde \Phi_{N} \pi_{+}\left( \frac{P}{g_{2}}\right) \right) \right).$$
  with
  \[
\left\{
\begin{array}{ccc}
  \Phi_{N}& =  &  \frac{g_{1}}{g_{2}} \chi^{N+1},\\
  \tilde \Phi_{N}& =  &  \frac{g_{2}}{g_{1}} \chi^{-(N+1)},  \\
  H_{\Phi_{N}} (\Psi) &  = &   \pi_{-} (\Phi_{N}\Psi) \\
   H^{*}_{\Phi_{N}} (\tilde \Psi) &  = &   \pi_{+} (\tilde \Phi_{N}\Psi)
\end{array}
\right.
\]
for  $\Psi \in \mathbb H^{+},$ (resp. $\tilde \Psi \in \left( \mathbb H^{+}\right)^{\bot}$)
and where $\pi_{+}$ (resp. $\pi_{-})$ are the orthogonal projection on 
$\mathbb H^{+}$ (resp. $\left(  \mathbb H^{+}\right)^{\bot}$).
\end{theoreme}
The reader can see \cite{RS020} for the statement and the proof of Theorem \ref{THEO2}.
In the appendix of this article we  briefly recall how to use it to calculate $\left( (T_{N})^{-1}(f_{0})\right)_{(1,1)}$ where the symbol $f_{0}$ is defined by 
$f_{0}=\chi_{0} \left(1-\bar \chi_{0}\chi \right)\left(1- \bar \chi_{0}\bar \chi  \right) 
\frac{1}{ \vert P_{N+1}\vert ^2}$. The equation 
(\ref{lastone}) gives the expression of 
$\left( (T_{N})^{-1}(f_{0})\right)_{(1,1)}$).
%%%%
We use this expression 
to obtain the equation (\ref {INV11}) which is a fundamental tool of our proof.\\
Always to obtain (\ref {INV11}) we have to use the fundamental property of the predictor polynomials
 which is the property (\ref{MAINPROP}).
  Before stating this property, we need of course to recall the definition of the predictor polynomial and its main property. 
\begin{definition}
The predictor polynomial of degree $M$ of a regular function $h$ is the trigonometric polynomial $K_M$ defined by 
$$ K_M= \sum_{k=0}^M \frac{ \left(T_{M}(h)\right)^{-1} _{k+1,1}}
{ \sqrt {\left(T_{M}(h) \right)^{-1}_{1,1}}} \chi ^k.$$
\end{definition}
\begin{prop}\label{MAINPROP}
For all integers $j$, such that $-M\le j \le M$ we have 
$$ \widehat{\left(\frac{1}{\vert K_M\vert ^2}\right)}(j)
= \hat h (j).$$ 
\end{prop}
We have also the useful  property 
\begin{prop}\label{SECONDPROP}
$K_{M} (r e^{i\theta}) \neq 0$ for all $\theta \in \mathbb R$ and $0\le r \le 1$.
\end{prop}
Finally, if $ \mathcal P_{M}$ is the set of trigonometric polynomials of degree less than or equal to $M$ we consider the scalar product defined on 
$ \mathcal P_{M}$ by 
$ \langle P \vert Q \rangle = \int _0^{2 \pi} P(\theta) \overline {Q(\theta)} h(\theta) d\theta$ and let's denote $\Phi_0, \Phi_1, \cdots, \Phi_M $ the orthogonal polynomials for this scalar product. 
Then the predictor polynomials of degree $0, 1, \cdots, M$ are closely related to these orthogonal polynomials by the relation 
$$ K_{j} (z) =  z ^j  \bar \Phi_{j}\left (\frac{1}{z}\right) \quad \forall j, 0 \le j\le M \quad \mathrm{and} \quad \forall z \neq 0.$$
The reader can consult \cite{Ld} for the predictor polynomials.\\
  We can now begin the demonstration of the theorem. \\
  This demonstration is divided into three parts. 
In the first part we obtain the equation (\ref{NEUF}) whose solutions are of the form $f^{-1}(\lambda)$ where the reals $\lambda$ are the eigenvalues of $f$ belonging to $[\theta_1, \theta_2]$.  In the second part we obtain an integral expression for the $\rho_N$ functions involved in this equation, which gives us the uniform convergence of the $\rho_N$ in $[\theta_1, \theta_2]$ to a continuous function $\rho$. The $\rho_N$
are therefore uniformly bounded, allowing us to locate the solutions of (\ref{NEUF}). In the third part, we study the smoothness of the function $\rho$, and use Taylor's theorem in (\ref{NEUF}) to obtain the asymptotic formula stated in Theorem.  
\subsection{Equation for the eigenvalues}
 Using the assumptions we can write $f(\theta) = f_{1}(1-\cos \theta)$
 where $f_{1}$ is a differentiable function strictly increasing on $[0,2]$. For all 
 $\lambda$ in $[f(\theta_{1}),f(\theta_{2})]$ we put 
 $\theta_{\lambda} = f^{-1} (\lambda)$ and $\lambda' = f_{1}^{-1}(\lambda)$, that means 
 $ \theta_{\lambda} = \arccos (1- \lambda')$, and 
 $\theta_{\lambda} \in [0, \pi]$.
 \begin{remarque}
 In the next proof we denote by $I_{\theta_{1},\theta_{2}}$ the set $[f(\theta_{1}), f(\theta_{2})]$.
 \end{remarque}
For $\lambda\in I_{\theta_{1},\theta_{2}}$ we have 
 $$ f(\theta) - \lambda=f_{1}(1-\cos \theta ) -\lambda=\left((1-\cos \theta)-(1-\cos \theta_{\lambda})\right) 
 H_{\lambda}(\theta)$$
 where $H_{\lambda} : \theta \mapsto H(\theta, \theta_{\lambda})$ is a regular function on $[-\pi, \pi]$.
   We can write 
 \begin{equation}\label{TROIS}
(1-\cos \theta)-(1-\cos \theta_{\lambda}) = (1-\cos \theta)-\lambda'= \frac{1}{2} \left(\vert 1-\chi\vert ^{2} - 2\lambda'\right).
 \end{equation}
  If $\chi_{\lambda} = e^{i \theta_{\lambda}}$ 
  we have  $\chi_{\lambda}= (1-\lambda') +i \sqrt { 1-(\lambda'-1)^{2}}$ and we can write the equation 
 (\ref{TROIS}) as 
 \begin{equation}\label{CINQ} 
(1-\cos \theta)-(1-\cos \theta_{\lambda})= -\frac{1}{2} \chi_{\lambda} 
 (1- \bar \chi_{\lambda}\chi)(1-\bar \chi_{\lambda} \bar \chi).
 \end{equation}
 Denote by $P_{N+1,\lambda}$ the predictor polynomial of degree $(N+1)$ $H_{\lambda}$. 
 The property (\ref{MAINPROP}) allows to write the equation  
\begin{equation} \label{EQUIVALENCE}
 T_{N} \Bigl( \left((1-\cos\theta)-(1-\cos\theta_{\lambda})\right)H_{\lambda'} \Bigr)=T_{N}\left( -\frac{1}{2} \chi_{\lambda} 
 (1- \bar \chi_{\lambda}\chi)(1-\bar \chi_{\lambda} \bar \chi) 
 \frac{1}{\vert P_{N+1,\lambda}\vert^{2}}\right).
 \end{equation} 
  For a fixed integer $N$ we denote by $T_{1,N,\lambda}$ the quantity  
  $ \left((T_{N}(f)-\lambda I_{N})^{-1}\right)_{1,1}$. \\
Since  $T_{N}(f) -\lambda I_{N}$ is a Toeplitz matrix we have, for  $\lambda \notin \operatorname{Spec}\left(T_{N}(f)\right)$
 $$T_{1,N,\lambda} = \frac{\det \left(T_{N-1}(f) -\lambda I_{N-1}\right)}
 {\det \left(T_{N}(f) -\lambda I_{N}\right)}.$$
 Using equation  (\ref{EQUIVALENCE}) and the inversion formula of Toeplitz matrices see in Theorem (\ref{THEO2}) we obtain the entry 
 $T_{1,N,\lambda}$.
  Then with the results (see the equation (\ref{lastone}) in the appendix) we can write, for $\lambda \notin  \operatorname{Spec}\left(T_{N}(f)\right)$,
  \begin{equation} \label{INV11}
 T_{1,N,\lambda} = \frac
 {\left(1-\bar \chi_{\lambda} ^{2(N+2)} \tau_{N}(\chi_{\lambda}) \right ) B_{2,N,\lambda} -B_{1,N,\lambda}}
 {1-\bar \chi_{\lambda} ^{2(N+1)} \tau_{N}(\chi_{\lambda})},
 \end{equation}
 and, for $\lambda \notin   \operatorname{Spec}\left(T_{N-1}(f)\right) \cup \operatorname{Spec}\left(T_{N}(f)\right)$ we have 
 \begin{equation}\label{NEWNEW}
 \frac{\det \left(T_{N}(f) -\lambda I_{N-1}\right)}
 {\det \left(T_{N-1}(f) -\lambda I_{N}\right)} = \frac  {1-\bar \chi_{\lambda} ^{2(N+1)} \tau_{N}(\chi_{\lambda})}
  {\left(1-\bar \chi_{\lambda} ^{2(N+2)} \tau_{N}(\chi_{\lambda}) \right ) B_{2,N,\lambda} -B_{1,N,\lambda}},
 \end{equation}
 with 
 $$ \tau_{N}(\theta_{\lambda}) = 
 \frac{\bar P_{N+1,\lambda} (\chi_{\lambda}) P_{N+1,\lambda} (\chi_{\lambda})}
 {\bar P_{N+1,\lambda} (\overline{\chi_{\lambda}}) P_{N+1,\lambda} (\overline{\chi_{\lambda}})},$$
 and 
 $$ B_{1,N,\lambda} = \Bigl \vert P_{N+1,\lambda} (0)  \Bigr \vert^{2} \tau_{N}(\theta_{\lambda}) { \bar \chi }^{2N+2}
 (1-\bar \chi_{\lambda}^{2}), B_{2,N,\lambda} =\frac{ \chi_{\lambda}}{\vert P_{N+1,\lambda}(0)\vert^2 }.$$
 Property \ref{SECONDPROP} implies that the quantities  $ B_{1,N,\lambda}$ and $ B_{2,N,\lambda}$ are well defined 
 and are different from zero for $\lambda \in I_{\theta_{1},\theta_{2}}.$\\
 Now the applications 
 $$ \lambda \mapsto  \frac{\det \left(T_{N}(f) -\lambda I_{N-1}\right)}
 {\det \left(T_{N-1}(f) -\lambda I_{N}\right)},  \lambda \mapsto \tau_{N}(\theta_{\lambda}),  \lambda \mapsto B_{1,N,\lambda}, 
   \lambda \mapsto B_{2,N,\lambda},$$
   are continuous on $I_{\theta_{1},\theta_{2}}\setminus  \operatorname{Spec}\left(T_{N-1}(f)\right)$. Hence we 
   can write the equality (\ref{NEWNEW}) 
   for all real $\lambda \in I_{\theta_{1},\theta_{2}} \setminus  \operatorname{Spec}\left(T_{N-1}(f)\right)$. 
    Since the eigenvalues of $\left(T_{N-1}(f) \right)$ are not in  $\operatorname{Spec} \left(T_{N}(f)\right)$ 
    (see \cite{RAIRO,BARBA}), we can write, according to Property \ref{SECONDPROP} :
 \begin{equation}\label{SIX}
 \lambda \in \left( \operatorname{Spec} \left(T_{N}(f)\right) \cap I_{\theta_{1},\theta_{2}}\right)\iff
  \chi_{\lambda}^{2(N+2)}= \tau_{N}(\theta_{\lambda}),\lambda \in I_{\theta_{1},\theta_{2}}.
  \end{equation}
  Since the function $H_{\lambda}$ is even, the constant 
$\tau_{N} (\theta_{\lambda})$ can be rewritten as
$$\tau_{N} (\theta_{\lambda}) = \left(\frac{P_{N+1,\lambda} ( \chi_{\lambda})}
{P_{N+1,\lambda} (\bar \chi_{\lambda})}\right)^2.
$$
As the function 
  $ \theta \mapsto \frac{P_{N+1,f(\theta)} (e^{-i\theta})}
{P_{N+1f(\theta)} (e^{i\theta})}$ is continuous from $[\theta_{1},\theta_{2}]$ to $\{z \vert \vert z\vert =1\}$ we have a function $\rho_{N}$ defined and continuous on $[\theta_{1}, \theta_{2}]$        
  such that 
  $\tau_{N}(\theta_{\lambda}) = e^{2 i \rho_{N}(\theta_{\lambda})}$.
  Then equation (\ref{SIX}) can be written 
    \begin{equation}\label{SEPT}
 \lambda \in \left(\mathcal \spec \left(T_{N}(f)\right) \cap[ f (\theta_{1}),f(\theta_{2})] \right) \iff
\theta_{\lambda} = \frac{\rho_{N} (\theta_{\lambda}) +  k \pi}{ (N+ 2)} , k \in[0, 2N+3].
  \end{equation}
  More precisely if $ M_{N} = \max_{\theta \in [\theta_{1}, \theta_{2}]}\vert \rho_{N} (\theta)\vert $ and if $\theta_{1}<a<b<\theta_{2}$
  we can write, according to the construction of $\chi_{\lambda}$    
   \begin{equation}\label{SEPT2}
 \lambda \in \left(\spec \left(T_{N}(f)\right) \cap [f(a), f(b)] \right) 
 \Rightarrow
 \theta_{\lambda} =\frac{ \rho_{N} (\theta_{\lambda}) +  k\pi}{N+ 2},
  k \in \Bigl[\frac{(N+2) a-M_{N}}{\pi}, \frac{(N+2) b+M_{N}}{\pi}\Bigr]
  \end{equation}
and 
 \begin{equation}\label{SEPTBIS}
 \theta_{\lambda} = \frac{\rho_{N} (\theta_{\lambda}) +  k\pi}{N+2}, 
  \frac{k \pi}{N+2}\in [a,b] \Rightarrow
 \lambda \in \left(\spec \left(T_{N}(f)\right) \cap \Bigl[f(a - \frac{M_{N}}{N+2}),
 f(b) + \frac{M_{N}}{N+2}\Bigr]\right).
 \end{equation}
  Lastly it is clear that we have now to solve the equation
   \begin{equation} \label{NEUF}
\theta = \frac{\rho_{N} (\theta) +  k\pi}{ (N+ 2)}.
\end{equation}

Now we have to make a more precise study of the function $\rho_{N}$. 
If $s\ge 0$, then every function $f \in \mathbb A(\mathbb T,s)$ without zeros on $\mathbb T$ admits a Wiener-Hopf factorization, that is, there exist functions $f_{+}$ and $f_{-}$
such that $f(e^{i\theta}) = f_{+}(e^{i\theta}) e^{i\gamma \theta}
f_{-}(e^{i\theta})$ with some $\gamma \in \mathbb Z$ the index of the factorization. The function $f_{+}$ (resp. $f_{-})$ 
belongs to set $\mathbb A(\mathbb T,s)_{+}$ (resp. $\mathbb A(\mathbb T,s)_{-}$)
where 
$$ \mathbb A(\mathbb T,s)_{+} = \lbrace f\in W^s \vert 
f(e^{i\theta}) = \sum_{j=0} ^{+\infty} \hat f (j) e^{i j \theta} \rbrace $$
and
$$  \mathbb A(\mathbb T,s)_{-} = \{  f\in W^s \vert 
f(e^{i\theta}) = \sum_{j=0} ^{+\infty} \hat f (-j) e^{-i j \theta} \}.$$ 
Here we have clearly 
$ \vert P_{N+1,\lambda}(e^{i\theta})\vert^2
=P_{N+1,\lambda}(e^{i\theta})
P_{N+1,\lambda}(e^{-i\theta})$
and 
$$\left( \frac{1}{\vert P_{N+1,\lambda}(e^{i\theta})\vert^2}\right)_{+}= \frac{1}{P_{N+1,\lambda}(e^{i\theta})}
\quad \mathrm{and} \quad  
\left( \frac{1}{\vert P_{N+1,\lambda}(e^{i\theta})\vert^2}\right)_{-}= \frac{1}{P_{N+1,\lambda}(e^{-i\theta})},$$
with index zero.
Now it is well known that in the Wiener-Hopf factorization 
$\left( \frac{1}{\vert P_{N+1,\lambda}(e^{i\theta})\vert^2}\right)_{+}$ can be written in the form 
$$ \left( \frac{1}{\vert P_{N+1,\lambda}(e^{i\theta})\vert^2}\right)_{+} = 
\exp \left ( \frac{1}{2} \log \Bigr(\frac{1}{\vert P_{N+1,\lambda}(e^{i\theta})\vert^2}\Bigl) + \frac{1}{2\pi i} P.V. \int _{\mathbb T}
 \frac{\log \left(\frac{1}{ \vert P_{N+1,\lambda}(z)\vert^2 } \right) }
 { z-e^{i\theta}} dz \right)$$
 that can be rewritten as 
\begin{align*}& \exp \left( \frac{1}{2} \log \Bigl(\frac{1}{\vert P_{N+1,\lambda}(e^{i\theta})\vert^2}\Bigr) + \frac{1}{4\pi i} P.V. \int _{0}^{2\pi} 
 \frac{\log \left(\frac{1}{ \vert P_{N+1,\lambda}(e^{i u}) \vert^2} \right)} 
 { \tan \frac{u-\theta}{2}} du \right. \\ &+
\left. \frac{1}{4\pi } \int _{0}^{2\pi} 
\log \left(\frac{1}{\vert P_{N+1,\lambda}(e^{i u}) \vert^2}\right)
du \right ).
\end{align*}
 That provides 
$ \frac{P_{N+1,\lambda}(e^{i\theta})} {P_{N+1,\lambda}(e^{-i\theta})} =  e^{i  \rho_{N,\lambda}(\theta)}$ with
$$  \rho_{N,\lambda}(\theta) = \frac{1}{4\pi }P.V.  \int _{0}^{2\pi} 
 \frac{\log \left(\frac{1}{ \vert P_{N+1,\lambda}(e^{i u})\vert^2 } \right)} 
 { \tan \frac{u-\theta}{2}} du - \frac{1}{4\pi} P.V.  \int _{0}^{2\pi} 
 \frac{\log \left(\frac{1}{ \vert P_{N+1,\lambda}(e^{i u}) \vert^2} \right)} 
 { \tan \frac{u+\theta}{2}} du $$
 and finally $ \rho_{N}(\theta_{\lambda})= \rho_{N,\lambda}(\theta_{\lambda})$, and 
 $ \rho_{N}(\theta)= \rho_{N, f(\theta)}(\theta)$.
 The same methods give, for $\mathrm{G_{\lambda} = \left( H_{\lambda}\right)_{+}}$
 $$ G_{\lambda} (e^{i\theta}) = \left( H_{\lambda}(\theta)\right)_{+} =
 \exp \left ( \frac{1}{2} \log \left(H_{\lambda} (v)\right) + \frac{1}{2\pi i} P.V. \int _{\mathbb T}
 \frac{\log \left(H_{\lambda} (v)\right) }
 { e^{iv}-e^{i\theta}} dv \right)$$
 and we obtain 
 $$
 \rho_{\lambda}(\theta)
=  \frac{1}{4\pi} P.V.\int _{0}^{2\pi} 
 \frac{\log \left(H_{\lambda}(u) \right)} 
 { \tan \frac{u-\theta}{2}} du -  \frac{1}{4\pi } P.V.\int _{0}^{2\pi} 
 \frac{\log \left(H_{\lambda}(u) \right)} 
 { \tan \frac{u+\theta}{2}} du,$$ 
 Then we put $$ \rho(\theta_{\lambda}) = \rho_{\lambda}(\theta_{\lambda}), 
\mathrm{and} \quad \rho(\theta) = \rho_{f(\theta)} (\theta).$$
  \subsection{Limit of the sequence  $(\rho_{N})_{(N\in \mathbb N)} $}
 Now we need  to relate the two functions $\rho_{N}$ and $\rho$
 and for this we have to obtain the following property 
 \begin{prop}\label{PROP}
 When $N$ goes to the infinity $\vert \rho_{N}(\theta) -\rho(\theta)\vert = O(\frac{1}{N^2})$
 uniformly in $\theta\in [\theta_{1},\theta_{2}]$.
 \end{prop}
 Lemmas 1 to 4 are devoted to the prove of this property.
 \begin{lemme}\label{UNO1}
  For all $\lambda \in I_{\theta_{1}},\theta_{2}$ the function 
   $\theta \mapsto H_{\lambda}(\theta)$    is in $C^{3}([0, 2\pi])$ and for all $j \in \{0,1,2,3\}$ 
   we have a real $K_{j}$ not depending on $\lambda$ 
   such that 
  $\Vert H_{\lambda}^{(j)}\Vert_{\infty}\le K_{j}.$
   \end{lemme}
   \begin{remarque}
   Here $\Vert h \Vert_{\infty}=\sup_{x\in [0,2 \pi]} \vert h(x)\vert$ for all bounded function $h$ on $[0,2\pi]$.
   \end{remarque}
   \begin{preuve}{of the lemma \ref{UNO1}}
If $t= 1-\cos \theta$ and $t_{\lambda}= 1- \cos \theta_{\lambda}$ we have to prove that the function 
$H_{1, \lambda} : t \mapsto \frac{f_{1}(t)- f_{1}(t_{\lambda})}{t-t_{\lambda}}$ is in $ C^{3} ([0, 2 \pi])$ and that for all 
integer $j$, $0\le j \le 3$ there exists a real $K_{1,j}$ such that,
for all $\lambda$ in $I_{\theta_{1},\theta_{2}}$ 
$\Vert H_{1, \lambda} ^{(j)}\Vert_{\infty}\le K_{1,j}.$  
Clearly $ \Vert H_{1,\lambda} \Vert_{\infty}\le\Vert  f^{(1)} \Vert_{\infty}$.
Now  for $t\neq t_{\lambda} $
\begin{align*}
H_{1,\lambda}^{(1)}(t) &= \frac{f_{1}^{(1)} (t) (t-t_{\lambda})
- \left( f_{1}(t) - f_{1}(t_{\lambda}) \right) } {(t-t_{\lambda})^2}\\
&= \frac{ \left (f_{1}^{(1)} (t_{\lambda})+ (t-t_{\lambda}) f_{1}^{(2)} (a_{1})\right) (t-t_{\lambda}) -\left( f_{1}^{(1)} (t_{\lambda}) 
(t-t_{\lambda}) + 
 \frac{ (t-t_{\lambda}) ^2}{2}  f_{1}^{(2)} (a_{2})\right)} {(t-t_{\lambda})^2}
 \end{align*}
 with $a_{1}$ and $a_{2}$ between $t$ and $t_{\lambda}$.
 That provides 
 \begin{itemize}
 \item
$ H_{1,\lambda}^{(1)}(t_{\lambda}) = \frac{ f_{1}^{(2)}(t_{\lambda})}{2},$
 \item
 $\Vert H_{1,\lambda}^{(1)}\Vert_{\infty}\le \frac{3}{2} 
 \Vert f_{1}^{(2)} \Vert_{\infty}$.
\end{itemize}
Now we have, for $t\neq t_{\lambda}$
\begin{align*}
H_{1,\lambda}^{(2)}(t) &= 
\frac{f_{1}^{(2)}(t) (t-t_{\lambda})^2 -2 \left( f_{1}^{(1)} (t) (t-t_{\lambda}) -\left( f_{1}(t) -f_{1}(t_{\lambda}) \right) \right) }
{(t-t_{\lambda})^3}\\
&=\frac{ \left( f_{1}^{(2)}(t_{\lambda})+ f_{1}^{(3)} 
(a_{3}) (t-t_{\lambda})\right) (t-t_{\lambda})^2 -
2 \left(d_{1,\lambda}(t) - d_{2,\lambda}(t)\right)}{(t-t_{\lambda})^3}
\end{align*}
with 
\begin{align*}
d_{1,\lambda}(t) &=  f_{1}^{(1)} (t_{\lambda}) (t-t_{\lambda}) 
+  f_{1}^{(2)}(t_{\lambda})(t-t_{\lambda})^2 +
f_{1}^{(3)} (a_{4}) \frac{(t-t_{\lambda})^3}{2}\\
d_{2,\lambda}(t) &=
 f_{1}^{(1)} (t_{\lambda}) (t-t_{\lambda}) + 
 f_{1}^{(2)}(t_{\lambda}) \frac{(t-t_{\lambda})^2} {2}
 + f_{1}^{(3)} 
(a_{5}) \frac{(t-t_{\lambda})^3}{6},
\end{align*}
and $a_{3},a_{4}, a_{5}$ between $t$ and $t_{\lambda}$. 
That provides 
\begin{itemize}
\item
$H_{1,\lambda}^{(2)} (t_{\lambda})= 
\frac{f^{(3)}(t_{\lambda})}{3}$,
\item
$\Vert H_{1,\lambda}^{(2)} \Vert_{\infty} 
\le\frac{7}{3} \Vert f_{1}^{(3)}\Vert_{\infty}.$
\end{itemize}
Finally we can write, for $t\neq t_{\lambda}$
\begin{align*}
H_{1,\lambda}^{(3)}(t) &=
\frac{ f_{1}^{(3)} (t) (t-t_{\lambda})^3 
-3 \left( f_{1}^{(2)} (t) (t-t_{\lambda})^2
-2 \left( f_{1}^{(1)} (t) (t-t_{\lambda}) -
\left( f_{1}(t) -f_{1}(t_{\lambda})\right) \right)\right)}
{(t-t_{\lambda})^4}\\
&= 
\frac{ \left(f_{1}^{(3)} (t_{\lambda}) + f_{1}^{(4)} (a_{6})(t-t\lambda)
\right) (t-t_{\lambda})^3 -3 \left( d_{3,\lambda} (t) 
-2  d_{4,\lambda}(t)\right)}{(t-t_{\lambda})^4}.
\end{align*}
with
\begin{align*}
d_{3,\lambda} (t) &= f_{1}^{(2)} (t_{\lambda}) (t-t_{\lambda})^2
+ f_{1}^{(3)} (t_{\lambda}) (t-t_{\lambda})^3 
+ f_{1}^{(4)}(a_{7}) \frac{(t-t_{\lambda})^4} {2} \\
d_{4,\lambda}(t) &= 
f_{1}^{(1)}(t_{\lambda)} (t-t_{\lambda}) + f_{1}^{(2)} (t_{\lambda}) (t-t_{\lambda})^2+
  f_{1}^{(3)} (t_{\lambda})\frac{ (t-t_{\lambda})^3 }{2}
  +  f_{1}^{(4)}(a_{8}) \frac{(t-t_{\lambda})^4} {6} 
\end{align*}
and $a_{6},a_{7}, a_{8}$ between $t$ and $t_{\lambda}$.
This last equalities give us 
\begin{itemize}
\item
$H_{1, \lambda}^{(3)}(t_{\lambda}) = \frac{1}{4} f_{1}^{(4)} (t_{\lambda}),$
\item
$\Vert H_{1,\lambda}^{(3)}\Vert _{\infty} \le 
\frac{15}{4} \Vert f_{1}\Vert_{\infty},$
\end{itemize}
which end the proof.
\end{preuve}
  \begin{remarque}\label{norme}
 If $h$ is a function is $L^2 ([0,2\pi])$   we denote by $\Vert h \Vert_{q,2} $ the quadratic norm $\left( \int_{0}^{2 \pi} \vert h((t) \vert ^2 dt \right)^{\frac{1}{2}}$.
 \end{remarque}
 \begin {lemme} \label{QUATRO}
 We have a real $S_{0}$ not depending on $k$ and $\lambda$
 such that 
 $$\Bigl \vert \widehat { G_{\lambda} }(k) \Bigr\vert \le \frac{S_{0}}{k^3},  \quad \Bigl \vert
\widehat { \frac{1}{G_{\lambda}} }(k) \Bigr\vert \le \frac{S_{0}}{k^3}, \quad \mathrm{for} \quad k>0, 
 \quad \mathrm{and} \quad  \Bigl\vert \widehat {\frac{G_{\lambda}}{\bar G_{\lambda}}}(k) \Bigr \vert \le\frac{S_{0}}{k^3},\quad \mathrm{for} \quad k\neq 0.$$
 \end{lemme}
 \begin{preuve}{}
 We can observe that, for $0\le j \le 3$, 
 $ \left(\pi_{+} \left (\log H_{\lambda}\right)\right)^{(j)}
 = \pi_{+}\left( \left (\log H_{\lambda}\right)^{(j)}\right)$.
 Hence, with Lemma \ref{UNO1}, we have, for  $0\le j \le 3$,
 \begin{equation}\label{GREAT}
 \Vert \left(\pi_{+} \left (\log H_{\lambda}\right)\right)^{(j)}
 \Vert_{q,2}\le \Vert \left (\log H_{\lambda}\right)^{(j)}\Vert_{q,2}
 \le T_{j}.
 \end{equation}
  If $m_{0} $ is the minimum of $H$ on $[0,2 \pi] \times [\theta_{1}, \theta_{2}]$ it is clear that for all $j \in \{0,1,2,3\}$
   $T_{j}$ is only depending on the constants $m_{0}, K_{0}, 
  K_{1}, K_{2}, K_{3},$ of Lemma \ref{UNO1}.
  Hence $T_{j}$ is no depending from $\lambda$.\\
 On the other hand since $\log H_{\lambda} \in C^3
 \left([0,2\pi]\right)$ we have, for all $n\ge 0$
 $$ \vert \widehat{ \log H_{\lambda}} (n)\vert \le 
\frac{ \Vert (\log H_{\lambda})^{(3)}\Vert_{q,2} }{n^3} \le
\frac{ T_{3}}{n^3}$$
and 
 $$ \vert \widehat{ \left(\log H_{\lambda} \right)^{(1)}}(n)\vert \le 
\frac{ \Vert (\log H_{\lambda})^{(3)}\Vert_{q,2} }{n^2} \le
\frac{T_{3}}{n^2}.$$
 Hence 
 \begin{equation} \label{GREAT2}
 \Vert\exp\left( \pi_{+}(\log H_{\lambda})\right)\Vert_{\infty}\le \exp \left( T_{3} \sum_{n\ge 0} \frac{1}{n^3}\right)=M_{1}
 \end{equation}
 and 
  $$ \Vert\left( \pi_{+}(\log H_{\lambda})\right)^{(1)} \Vert_{\infty}\le  T_{3} \sum_{n\ge 0} \frac{1}{n^2}= M_{2}.$$
 Now if we put 
 $ \pi_{+}(\log H_{\lambda}) = F_{\lambda}$
 we can write 
 $$ \left(\exp -F_{\lambda}\right)^{(3)}=
 \left(- F_{\lambda}^{(3)}  + 3 F_{\lambda}^{(1)} 
 F_{\lambda}^{(2)} - \left(F_{\lambda}^{(1)}\right)^3\right) 
 \exp -F_{\lambda},$$
 and 
 $$ \left(\exp F_{\lambda}\right)^{(3)}=
 \left( F_{\lambda}^{(3)}  + 3 F_{\lambda}^{(1)} 
 F_{\lambda}^{(2)}+ \left(F_{\lambda}^{(1)}\right)^3\right) 
 \exp F_{\lambda}.$$
 According to (\ref{GREAT}) we have the inequalities, 
 $$ \Vert F_{\lambda}^{(3)}\exp F_{\lambda}\Vert_{q,2}
 \le \Vert F_{\lambda}^{(3)}\Vert_{q,2} \Vert \exp F_{\lambda}\Vert_{\infty}\le T_{2}M_{1},$$
 $$ \Vert F_{\lambda}^{(1)} F_{\lambda}^{(2)}\exp F_{\lambda}\Vert_{q,2} \le \Vert F_{\lambda}^{(1)}\Vert _{\infty}
\Vert \exp F_{\lambda}\Vert_{\infty}\Vert F_{\lambda}^{(2)}\Vert _{q,2}\le M_{1}M_{2}T_{2},$$
$$ \Vert (F_{\lambda}^{(1)})^3\exp F_{\lambda}\Vert_{q, 2}
\le \Vert (F_{\lambda}^{(1)})\Vert_{\infty}^3
\Vert \exp F_{\lambda}\Vert_{\infty}\le M_{2}^3 M_{1}.$$
 This means that 
$\Vert \left(\exp F_{\lambda}\right)^{(3)}\Vert_{2}$,  
is bounded by a constant $S_{1}$ not depending on $\lambda$ and $n$.
This result implies 
$$ \vert \widehat { \exp- F_{\lambda}} (n) \vert \le 
\frac{ \Vert \left(\exp -F_{\lambda}\right)^{(3)}\Vert_{q,2}}{n^3}
\le \frac{S_{1}}{n^3}$$
and 
$$ \vert \widehat { \exp F_{\lambda}} (n) \vert \le 
\frac{ \Vert \left(\exp F_{\lambda}\right)^{(3)}\Vert_{q,2}}{n^3}
\le \frac{S_{1}}{n^3}$$
for all $n\ge0$, that is the first part of the lemma.\\
On the other hand for $n > 0$ we have 
\begin{align*}
\Bigl \vert \widehat{ \frac{ G_{\lambda}}{\bar G_{\lambda}}} (n)\Bigr \vert  &= 
\Bigl \vert \sum_{h\ge 0} \widehat {G_{\lambda}} (h+n) \widehat {\frac{1}{\bar G_{\lambda}}}(-h) \Bigr \vert \\
& \le S_{1}^2 \frac{ \sum _{h>0} \frac{1}{h^3}} {n^3} + \frac{ 1}{n^3} \vert\widehat{ \frac{1}{\bar G_{\lambda}}}(0)\vert,
\end{align*}
and 
\begin{align*}
\Bigl \vert \widehat{ \frac{ G_{\lambda}}{\bar G_{\lambda}}} (-n)\Bigr \vert  &= 
\Bigl \vert \sum_{k\ge 0} \widehat {G_{\lambda}} (k) \widehat {\frac{1}{\bar G_{\lambda}}}\left(-(k+n)\right) \Bigr \vert \\
& \le S_{1}^2 \frac{ \sum _{k>0} \frac{1}{k^3}} {n^3} + \frac{ 1}{n^3} \vert\widehat{ \frac{1}{\bar G_{\lambda}}}(0)\vert,
\end{align*}
and with (\ref{GREAT2}) we can write 
$\vert\widehat{ \frac{1}{\bar G_{\lambda}}}(0)\vert \le 
\frac{1}{2 \pi} \Vert \exp ( \pi_{+} ( \log H_{\lambda}) \Vert_\infty\le \frac{M_{1}}{2 \pi}
$
that provides the third inequality of the lemma. 
\end{preuve} \begin{lemme} \label{CINQUE}
 If $ \beta_{k+1,\lambda} =  \widehat { \frac{1}{G_{\lambda}} }(k)$ we have, for a sufficient large $N$
 $$(\left( T_{N} \left( H_{\lambda}\right) \right)^{-1} _{k,1} = \overline{ \beta_{0,\lambda}} 
 \beta_{k,\lambda}+R_{k,N,\lambda}$$
 with 
 $ \vert R_{k,N,\lambda} \vert \le \frac{M}{N^2(N+1-k)^2}$ where $M$ is not depending on  $\lambda$ and $k$.
 \end{lemme}
 \begin{preuve}{}
 Using the inversion formula given in the appendix of this paper we obtain, for \mbox{$H_{\lambda}= G_{\lambda} \bar G_{\lambda}$,}
$G_{\lambda} \in \mathbb H^{+}$,
$$ \left( T_{N}(H_{\lambda})\right)^{-1}_{k+1,l+1} = \Big\langle \pi_{+}\left( \frac{\chi^{l}}{\bar G_{\lambda} }\right)\vert \frac{\chi^{k}}{\bar G_{\lambda} } 
\Big\rangle - \Big\langle \sum_{s=0}^{+\infty} \left( H^{*}_{\Phi_{N,\lambda}} H_{\Phi_{N,\lambda}}\right)^{s} \pi_{+} \bar \Phi_{N,\lambda}
\pi_{+} \left( \frac{\chi^{l}}{\bar G_{\lambda} } \right )\vert  \pi_{+} \bar \Phi_{N,\lambda}
\pi_{+} \left( \frac{\chi^{k}}{\bar G_{\lambda} }\right)\Big\rangle,$$
with $$\Phi_{N,\lambda} = \frac {G_{\lambda}}{\bar G_{\lambda}} \chi ^{N+1}, \quad \mathrm{and} \quad
\tilde \Phi_{N,\lambda} = \frac{\bar G_{\lambda}}{ G_{\lambda}} \chi ^{-(N+1)},$$
$$ H_{\Phi_{N,\lambda}} (\Psi)=\pi_{-}(\Phi_{N,\lambda} \Psi)
\quad \mathrm{for} \quad  \Psi \in \mathbb H^+, $$
$$ H^{*}_{\Phi_{N,\lambda}} (\Psi)=\pi_{+}(\tilde \Phi_{N,\lambda} \Psi) \quad 
\mathrm{for} \quad  \Psi \in(\mathbb H^{+})^{\perp}.$$
 For $l=0$ this formula becomes 
$$ \left( T_{N}(H_{\lambda})\right)^{-1}_{k+1,1}= \Big\langle \pi_{+}\left( \frac{1}{\bar G_{\lambda} }\right)\vert \frac{\chi^{k}}{\bar G_{\lambda} } 
\Big\rangle -\Big \langle \sum_{s=0}^{+\infty} \left( H^{*}_{\Phi_{N,\lambda}} H_{\Phi_{N,\lambda}}\right)^{s} \pi_{+} 
\bar \Phi_{N,\lambda}
\pi_{+} \left( \frac{1}{\bar G_{\lambda} }\right )\vert  \pi_{+} \bar \Phi_{N,\lambda}
\pi_{+} \left( \frac{\chi^{k}}{\bar G_{\lambda} }\right)\Big\rangle.$$
In the next of this proof we use the following notation : 
$$\frac{G_{\lambda}}{\bar G_{\lambda}} = \sum_{u \in \mathbb Z} \gamma_{u,\lambda }\chi^u .$$
From Lemma \ref{QUATRO} we have a positive constant  $S_{0}$  such that 
$$ \vert \beta_{u,\lambda}\vert \le \frac{S_{0}}{u^3} \quad \forall u\in \mathbb N^\star\quad 
\mathrm{and} \quad \vert \gamma_{u,\lambda}\vert \le \frac{S_{0}}{u^3} \quad \forall u\in \mathbb Z^\star.$$
We obtain 
$$\Big \langle \pi_{+}\left( \frac{1}{\bar G_{\lambda} }\right)\vert \frac{\chi^{k}}{\bar G_{\lambda} }\Big \rangle =
\bar \beta_{0,\lambda} \beta_{k,\lambda},$$
$$ \pi_{+} \bar \Phi_{N,\lambda}
\pi_{+} \left( \frac{1}{\bar G_{\lambda} } \right ) = \pi_{+} \left( \bar \Phi_{N,\lambda} \bar \beta_{0,\lambda}\right) = 
 \bar \beta_{0,\lambda} \sum_{v\ge N+1 } \bar \gamma_{-v,\lambda} \chi^{v-N-1},$$
$$ \pi_{+} \bar \Phi_{N,\lambda}
\pi_{+} \left( \frac{\chi^{k}}{\bar G_{\lambda} }\right)= \sum_{w=0}^k \bar \beta_{w,\lambda} \left( \sum_{v\ge N+1-k+w}
\bar \gamma_{-v,\lambda}\chi^{v-N-1+k-w}\right).$$
Hence
we obtain
$$ \Bigl\Vert \pi_{+} \bar \Phi_{N,\lambda}
\pi_{+} \left( \frac{1}{\bar G_{\lambda} } \right )\Bigr\Vert _{q,2}
\le S_{1}\left (N+1\right)^{-2}, $$
and 
$$
\Bigl\Vert \pi_{+} \bar \Phi_{N,\lambda}
\pi_{+} \left( \frac{\chi^{k}}{\bar G_{\lambda} }\right)\Bigr \Vert_{q,2}
\le  S_{1}(\left (N+1-k\right)^{-2}, $$
where $S_{1}$ no depending on $\lambda$ and $N$
On the other hand for $\psi = \displaystyle{ \sum_{w\ge 0} \alpha_{w}\chi^w }$ a function in $ \mathbb H^+$ we have, with the continuity of the projection $\pi_{-}$, 
$$ H_{\Phi_{N,\lambda}} (\psi) = \sum_{w\ge 0} \alpha_{w} \left( \sum_{v>N+1+w} \gamma_{-v,\lambda} \chi^{-v+w+N+1}\right)$$
that provides
\begin{align*}
\Vert H_{\Phi_{N,\lambda}} (\psi) \Vert _{q,2} &\le \sum_{w\ge 0} \vert \alpha_{w}\vert \left (\sum_{v>N+1+w } \vert \gamma_{-v,\lambda} \vert 
\right)\\ 
& \le \Vert \psi \Vert _{2} \left(\sum_{w\ge 0} \left(\sum_{v>N+1+w } \vert \gamma_{-v,\lambda} \vert\right)^2\right)^{1/2}\\
&\le S_{0} \Vert \psi \Vert _{2} (N+1)^{-3/2}
\end{align*}
 Clearly we have also 
$\Vert H^\star_{\Phi_{N,\lambda}} (\psi)\Vert_{q,2} \le S_{0} (N+1)^{-3/2} \Vert (\psi)\Vert_{q,2}$ and we can write 
$$\Bigr \Vert \sum_{s=0}^{+\infty} \left( H^{*}_{\Phi_{N,\lambda}} H_{\Phi_{N,\lambda}}\right)^{s} \pi_{+} \bar \Phi_{N,\lambda}
\pi_{+} \left( \frac{1}{\bar G_{\lambda} }\right )\Bigl\Vert_{q,2} \le 
\frac{S_{1}} {\left(1- S_{0}^2 (N+1)^{-3}\right)^2} (N+1)^{-2}.$$
And finally 
we can write 
$$ \left( T_{N}(H_{\lambda})\right)_{1,k+1} ^{-1}= \bar \beta_{0,\lambda} \beta_{k,\lambda} 
+ O \left((N+1)^{-2}(N+1 -k)^{-2}\right)$$
with, for a sufficiently large $N$, $O \left((N+1)^{-2}(N+1 -k)^{-2}\right) = 2 S_{1}^2 
(N+1)^{-2}(N+1 -k)^{-2}$ 
uniformly in $\lambda$ that is the expected result.
\end{preuve}
\begin{remarque}
As the coefficient $\beta_{0,\lambda}$ is real the form 
of $\tau_{N}(\chi_{\lambda}) $ allows to assume that 
$\beta_{0,\lambda}=1$ is the rest of our demonstration.
\end{remarque}

 \begin{lemme} \label{septimus}
 We have 
 $ \Bigl\Vert \ln \left(\frac{1}{\vert P_{N+1,\lambda}\vert^2}\right)
  - \ln \left( H_{\lambda}\right)\Bigr\Vert_{0} = O\left( \frac{1}{N^2}\right) $ uniformly in $\lambda$.
  \end{lemme}
  \begin{preuve}{} 
   Using Lemma 
 \ref{CINQUE}, we obtain 
  $$\Bigl \Vert P_{N+1,\lambda} - \frac{1}{G_{\lambda} }\Bigr\Vert _{0}
  \le M (N+1)^{-2}\sum_{k=0} ^N\frac{1}{(N+1-k)^{2}}
  + \sum_{k=N+1} ^{+ \infty} \vert \beta_{k,\lambda} \vert.$$
 Hence 
   \begin{equation} \label{DOUDOU} \Bigl\Vert P_{N+1,\lambda} - \frac{1}{G_{\lambda} }\Bigr\Vert _{0}
  \le \frac{M+S_{0}}{(N+1)^2},
  \end{equation}
  where $M$ and $\beta_{k,\lambda}$ as in Lemma \ref{CINQUE}
  and $S_{0}$ is the real not depending on$N$ and from $\lambda$ which has been introduced in Lemma \ref{QUATRO}.
   Always with $M$ and $S_{0}$ no depending from $\lambda$  and  the norm $\Bigl \Vert P_{N+1,\lambda} - \frac{1}{G_{\lambda} }\Bigr\Vert _{0}$
   is bounded by $O\left( \frac{1}{N^2}\right)$.
   Now since $\Vert \Psi \Phi\Vert _{0} \le \Vert \Psi \Vert _{0} \Vert  \Phi\Vert _{0}$
we have 
\begin{equation} \label{DOUDOU1}
 \Bigl\Vert \frac{1}{P_{N+1,\lambda}} - G_{\lambda}\Bigr\Vert_{0} \le 
\Bigl\Vert P_{N+1,\lambda} - \frac{1}{G_{\lambda} }\Bigr\Vert _{0} 
\Bigl\Vert \frac{G_{\lambda}}{P_{N+1,\lambda}}\Bigr\Vert_{s_{0}}\le 
\Bigl\Vert P_{N+1,\lambda} - \frac{1}{G_{\lambda} }\Bigr\Vert _{0} \Bigl\Vert \frac{1}{P_{N+1,\lambda}}\Bigr\Vert_{0}
\Bigl\Vert G_{\lambda}\Bigr\Vert_{0}.
\end{equation}
Then, according to (\ref{DOUDOU1}) we have
\begin{equation}  \label{DOUDOU2}  
\Bigl\Vert \frac{1}{P_{N+1,\lambda}} \Bigl\Vert_{0} -  \Bigl\Vert G_{\lambda}\Bigr\Vert_{0}
\le \Bigl\Vert P_{N+1,\lambda} - \frac{1}{G_{\lambda} }\Bigr\Vert _{0} 
\Bigl\Vert \frac{1}{P_{N+1,\lambda}}\Bigr\Vert_{0}
\Bigl\Vert G_{\lambda}\Bigr\Vert_{0}
\end{equation}
That provides 
\begin{equation}  \label{DOUDOU3}  
\Bigl\Vert \frac{1}{P_{N+1,\lambda}}\Bigr\Vert_{0}
\left (1-\Bigl\Vert P_{N+1,\lambda} - \frac{1}{G_{\lambda} }\Bigr\Vert _{0}
\Bigl\Vert G_{\lambda}\Bigr\Vert_{0} \right) 
\le \Vert G_{\lambda}\Vert_{0}
\end{equation}
According to Lemma  \ref{QUATRO} we have a real $A_{1}$ such that for all $\lambda$ in $]f(\theta_{1}, f(\theta_{2}[$ we have
$ \Vert G_{\lambda}\Vert_{0}
 \le A_{1}$. Hence with (\ref{DOUDOU}) we obtain that for $N$ sufficiently large we have 
 \begin{equation} \label{DOUDOU4}
1-\Bigl\Vert P_{N+1,\lambda} - \frac{1}{G_{\lambda} }\Bigr\Vert _{0}
\Bigl\Vert G_{\lambda}\Bigr\Vert_{0}
\ge \frac{1}{2}.
\end{equation}
and 
\begin{equation} \label{DOUDOU5}
\Bigl\Vert \frac{1}{P_{N+1,\lambda}}\Bigl\Vert_{0} \le 2A_{1}.
\end{equation}
Merging (\ref{DOUDOU}) and (\ref{DOUDOU5} )we obtain 
\begin{equation}\label{DOUDOU6}
 \Bigl\Vert \frac{1}{\vert P_{N+1,\lambda}\vert^2} - H_{\lambda}\Bigr\Vert _{0}
\le\Bigl \Vert \frac{1}{\vert P_{N+1,\lambda}\vert^2} - \frac{1}{P_{N+1,\lambda}} \bar G_{\lambda}\Bigr\Vert _{0}
+\Bigl \Vert  \frac{1}{P_{N+1,\lambda}} \bar G_{\lambda} -H_{\lambda}\Bigr\Vert_{0}\le 
4 A_{1}^3 \left( \frac{M+S_{0}}{(N+1)^2}\right)
\end{equation}
that means
 $\Bigr\Vert \frac{1}{\vert P_{N+1,\lambda}\vert^2} -H_{\lambda}\Bigl\Vert_{0}=O\left( \frac{1}{N^2}\right)$,
 uniformly in $\lambda$. Now observe that 
   $$
 \Bigl  \Vert \ln \left(\frac{1}{\vert P_{N+1,\lambda}\vert^2}\right)
 - \ln \left(  H_{\lambda}\right)\Bigr\Vert_{0} = 
 \Bigl \Vert \ln \left( 1+\frac{\frac{1}{\vert P_{N+1,\lambda}\vert^2} -H_{\lambda}}{H_{\lambda}}\right)\Bigr\Vert_{0},
  $$
  that is also 
  $$ 
  \Bigl \Vert \ln \left( 1+\frac{\frac{1}{\vert P_{N+1,\lambda}\vert^2} -H_{\lambda}}{H_{\lambda}}\right)\Bigr\Vert_{0}
  \le  \sum_{n\ge1} \frac{1}{n} \left(\Bigl\Vert \frac{1}{\vert P_{N+1,\lambda}\vert^2} -H_{\lambda}\Bigr\Vert_{0}\right)^n
\left(  \Bigl\Vert \frac{1} {H_{\lambda}}\Bigr\Vert _{0}\right)^n.$$
  Now we have, according to Lemma \ref{UNO1}, 
  $$ \Bigr\Vert \frac{1} {H_{\lambda}}\Bigl\Vert _{0}
    \le \left(\sum_{n\ge 0} \frac{1}{n^2} \right) \frac{1}{2\pi}\Bigl\Vert \left(\frac{1}{H_{\lambda}}\right)^{(2)}\Bigr \Vert _{q,2}
  \le \frac{K}{m_{0}^3}
  $$
  with 
  $m_{0}$ as in the proof of Lemma \ref{QUATRO} and $K$ no depending on $\lambda$ and $N$. 
  That gives us, according to (\ref{DOUDOU6})
 $$ \Bigl \Vert \ln \left( 1+\frac{\frac{1}{\vert P_{N+1,\lambda}\vert^2} -H_{\lambda}}{H_{\lambda}}\right)\Bigr\Vert_{0}
  \le  \sum_{n\ge1} \frac{1}{n} \left(\Bigr\Vert 4 A_{1}^3 \left( \frac{M+S_{0}}{(N+1)^2}\right)
\Bigl\Vert_{0}\right)^n
\left( \frac{K} {m_{0}^3}\right )^{m}.$$
 Since $m_{0}$ and $K$ are not depending on  $\lambda$ we can conclude 
$$ \Bigl  \Vert \ln \left(\frac{1}{\vert P_{N+1,\lambda}\vert^2}\right)
 - \ln \left(  H_{\lambda}\right)\Bigr\Vert_{0} = 
 \Bigl \Vert \ln \left( 1+\frac{\frac{1}{\vert P_{N+1,\lambda}\vert^2} -H_{\lambda}}{H_{\lambda}}\right)\Bigr\Vert_{0}=O\left(\frac{1}{N^2}\right)
  $$
uniformly in $\lambda$.
  \end{preuve}
 Since the Cauchy  singular operator is bounded on the Wiener classes $A(\mathbb T,s)$, $s\ge0$,
 we have $\Vert \rho_{N}  - \rho \Vert _{0} = O(\frac{1}{N^2})$ and $\vert \rho_{N}(\lambda) -\rho(\lambda)\vert = O(\frac{1}{N^2})$
 uniformly in $\lambda$. That ends the proof of Property \ref{PROP}.
 \subsection{Derivation and solutions of the equation for the eigenvalues}
To do this we need the two following lemmas.
\begin{lemme} \label {L1}
The function $\rho$ is in $\mathcal C^2 ([\theta_{1},\theta_{2}])$.
\end{lemme}
\begin{preuve}{}
We prove the result for the function
$$ I : \theta \mapsto P.V. \int_{0}^{2 \pi}
 \frac{\ln \left( H(t,\theta)\right)}{ \tan \left( \frac{t-\theta}{2}\right)} 
 dt,$$
 the proof is quite the same for the function 
 $$ \theta \mapsto P.V. \int_{0}^{2 \pi}
 \frac{\ln \left( H(t,\theta)\right)}{ \tan \left( \frac{t+\theta}{2}\right) }
 dt.$$
First we write 
$I(\theta) = I_{1,\theta} + I_{2,\theta}$ with
$$ I_{1,\theta} = P.V. \int_{0}^{2\pi}  \frac{\ln \left( H(\theta,\theta)\right)}{ \tan \left( \frac{t-\theta}{2}\right)} 
 dt,$$
 $$ I_{2,\theta}=  \int_{0}^{2\pi} 
 \frac{\log \left(H(t,\theta)\right)- \log \left(H(\theta,\theta)\right)} { \tan \left( \frac{t-\theta}{2}\right)} dt.$$
A simple calculus provides us $I_{1,\theta}=0$. On the other hand we can observe that 
the function $\Psi$ defined by 
$$\Psi(t,\theta) = \frac{\log \left(H(t,\theta)\right)- \log \left(H(\theta,\theta)\right)} { \tan \left( \frac{t-\theta}{2}\right)}$$
 can be write 
$$\Psi(t,\theta) =\frac{\log \left(H(t,\theta)\right)- \log \left(H(\theta,\theta)\right)} {\frac{t-\theta}{2}}
\frac{ {\frac{t-\theta}{2}}} { \tan \left( \frac{t-\theta}{2}\right)}.$$
Thanks to the symmetry of the function $H: \theta\mapsto H(\theta,\theta')$ we can say that the function 
$\theta \mapsto H (t,\theta)$ is in $\mathcal C^3\left ( [\theta_{1},\theta_{2}]\right)$ for all $t$ in 
$[0,2 \pi]$. Hence if $\Psi_{1}$ is the function defined by 
$\Psi_{1} : \theta \mapsto \frac{\log \left(H(t,\theta)\right)- \log \left(H(\theta,\theta)\right)} {t-\theta}$
the function $\frac{\partial \Psi_{1}}{\partial \theta} (t,\theta)$ is defined for all $\theta\neq t$ 
and is equal to 
$$ \frac{\left (t-\theta) \left( (\log H)'_{\theta} (t,\theta) -(\log H)'_{t} (\theta,\theta) -(\log H)'_{\theta}(\theta,\theta)
\right) + \left( (\log H)(t,\theta) -(\log H)(\theta,\theta)\right) \right)} {(t-\theta)^2},$$
where we have denoted by $(\log H)'_{t}$ the quantity $\frac{ \partial (\log H)}{\partial t}$ and by 
$(\log H)'_{\theta}$ the quantity $\frac{ \partial (\log H)}{\partial \theta}$
We see that for $t=\theta$ the function $\frac{\partial \Psi_{1}}{\partial \theta}$ is equal to 
$\frac{\partial^2\log H}{\partial t^2}(\theta,\theta)$  .
Since the functions $ \log H$, $ \frac{ \partial (\log H)}{\partial t}$, $ \frac{ \partial (\log H)}{\partial \theta}$, 
and $\frac{\partial^2\log H}{\partial \theta^2}$ are continuous on $[0,2 \pi] \times [\theta_{1},\theta_{2}]$
we obtain that the function $ (t,\theta)\mapsto \frac{\partial \Psi}{\partial \theta} $ is defined and continuous 
on $[0,2 \pi] \times [\theta_{1}, \theta_{2}]$, that completes this demonstration for the existence of $\rho^{(1)}$. For 
$\rho^{(2)}$ the function $ \frac{\partial^2 \Psi_{1}} {\partial \theta^2} (t,\theta) $ is defined for all $t\neq \theta$ and is 
equal to $$ \frac{ (t-\theta)^2  \Psi_{1}(t,\theta) +2 (t-\theta) \Psi_{3} (t, \theta) +2 \Psi_{4}(t, \theta)} {(t-\theta)^3}$$
where, 
\begin{align*}
 \Psi_{2}(t,\theta) &= (\log H)^{\prime \prime} _{\theta^2} ( t,\theta) - (\log H )^{\prime\prime}_{t^2} (\theta,\theta) 
 - (\log H )^{\prime\prime}_{\theta^2} (\theta,\theta)  -  (\log H )^{\prime\prime}_{t,\theta} (\theta,\theta), \\
\Psi_{3}(t,\theta) &= (\log H)'_{\theta} (t,\theta) -(\log H)'_{t} (\theta,\theta) -(\log H)'_{\theta}(\theta,\theta),\\
\Psi_{4}(t,\theta) &= (\log H)(t,\theta) -(\log H)(\theta,\theta),
\end{align*}
and for $t=\theta$ we see that $$\frac{\partial^2 \Psi_{1}}{\partial \theta^2} =  (\log H)^{(3)}_{\theta^2,t}
(\theta,\theta) +  (\log H)^{(3)}_{t^3} (\theta,\theta) + (\log H)^{(3)}_{t^2\theta}(\theta,\theta).$$
Then the same arguments as previously allow us to conclude.
\end{preuve}
To begin stating Theorem \ref{THEO1ter} we have to remark that with Property \ref{PROP} we have a real $M>0$ such that 
$ -M\le \rho_{N}(\theta) \le M$ for all integer $N$ and all $\theta \in [\theta_{1}, \theta_{2}]$.\\
Now if $\lambda\in [f(a), f(b)]$ is an eigenvalue of $T_{N}(f)$ we know that there is a real $\theta_{\lambda} \in [a,b]$ 
such that $\theta_{\lambda}$ is a solution of (\ref{NEUF}) that implies $ a-\frac{M}{N+2} \le \frac {k \pi}{N+2}\le b+\frac{M}{N+2}$, 
and we can conclude $k \in [k_{\theta_{1,N}},k_{\theta_{2,N}}]$ for $N$ sufficiently large.\\
 Reciprocally if $N$ is sufficiently large we have for all $k \in [k_{a,N},k_{b,N}] $  two reals $\theta'_{k}$ and $\theta^{\prime\prime}_{k}$ in 
$[\theta_{1}, \theta_{2}]$ such that 
$ \theta_{k}' < \frac{k\pi-M}{N+2} $ and $\frac{k\pi+M}{N+2}< \theta^{\prime\prime}_{k}$ that provides a solution to the equation 
$\theta = \frac{\rho_{N}(\theta) +k \pi}{N+2}$.\\
Now we can  obtain the formula announced in the statement of  Theorem \ref{THEO1ter}. %%%%\`u
 For 
$\lambda$ an eigenvalue in $]f (\theta_{1}), f(\theta_{2})[$ 
we have, 
following the equation (\ref{NEUF}), 
$ \lambda = f\left(\frac{k\pi+\rho_{N} (\theta_{\lambda}) }{N+2}\right)$
where $\theta_{\lambda}$ is a  solution of
 the equation (\ref{NEUF}).
According to Property \ref{PROP} we can enunciate 
\begin{prop} \label{beautiful}
$ \lambda = f\left(\frac{k\pi+\rho (\theta_{\lambda}) }{N+2}+\frac{R_{N,\lambda}}{(N+2)}\right)$
 with $R_{N,\lambda}= O \left( \frac{1}{(N+2)^2}\right)$ uniformly in 
 $\lambda$. 
 \end{prop}
Putting $d= \frac{\pi k}{N+2}$ we have by Taylor's theorem,
\begin{align*} 
 \lambda &= f(d) + f'(d) \left( \frac{\rho (\theta_{\lambda})+R_{N,\lambda} }{N+2}\right) +
  \frac{1}{2}
f^{\prime\prime} (d)  \left( \frac{\rho (\theta_{\lambda})+ R_{N,\lambda}}{N+2}\right)^2 \\
&+\frac{1}{6} f^{(3)} \left (d +h_{1} \frac{\rho (\theta_{\lambda}) +R_{N,\lambda}}{N+2}\right) 
\left(\frac{ \rho(\theta_{\lambda})+R_{N,\lambda}}{N+2}\right)^3,
\end{align*}
with $ 0<h_{1}<1$. 
That provides 
\begin{equation} \label{DIX}
 \lambda = f(d) + f'(d) \left( \frac{\rho (\theta_{\lambda}) }{N+2}\right) +
  \frac{1}{2}
f^{\prime\prime} (d)  \left( \frac{\rho (\theta_{\lambda})}{N+2}\right)^2 
+\frac{1}{6} f^{(3)} \left (d +h_{1} \frac{\rho (\theta_{\lambda}) }{N+2}\right) 
\left(\frac{ \rho(\theta_{\lambda})}{N+2}\right)^3 + O\left(\frac{1}{(N+2)^3}\right),
\end{equation}
where the quantity $ O\left( \frac{1}{(N+2)^3}\right)$ is bounded uniformly in $\lambda$.
On the other hand, with the equation (\ref{NEUF}) $\theta_{\lambda} = \frac{k\pi+\rho (\theta_{\lambda}) }{N+2}
+\frac{R_{N,\lambda}}{(N+2)^2}$
and we can write, always by Taylor's theorem, 
\begin{equation}\label{ONZE}
\rho(\theta_{\lambda})=  \rho (d) + \rho'(d)  \frac{\rho (\theta_{\lambda}) +R_{N,\lambda}}{N+2} +
\frac{1}{2} \rho^{\prime \prime} \left((d +h_{2} \frac{\rho (\theta_{\lambda})+R_{N,\lambda} }{N+2}\right) 
\left(\frac{ \rho(\theta_{\lambda})+R_{N,\lambda}}{N+2}\right)^2, \end{equation}
with $0<h_{2}<1$,
that implies 
\begin{equation}\label{ONZEBIS}
\rho(\theta_{\lambda})=  \rho (d) + \rho'(d)  \frac{\rho (d)}{N+2} + O\left(\frac{1}{(N+2)^2}\right).
\end{equation}
with the rest is bounded by $\frac{\vert S\vert }{(N+2)^2}$ where $S$ is a constant no depending from $\lambda$.
Merging the equation (\ref{DIX}) and  (\ref{ONZEBIS})
we obtain
\begin{equation} \label{DOUZE}
 \lambda =  f(d) + \frac{f'(d) \rho (d)}{N+2}  + \frac{ f'(d) 
\rho'(d) \rho(d) }{(N+2)^2} + \frac{1}{2}
\frac{ f^{\prime\prime} (d)  \rho^2 (d) }{(N+2)^2} 
+R_{N,d},
\end{equation}
with $ R_{N,d} = O \left( \frac{1}{(N+2)^3}\right)$ uniformly in $\lambda$.
To achieve the proof  we have to be sure that the eigenvalues found are distincts as announced. To do this we need the  following  three lemmas.
\begin{lemme} \label{L2}
For $k,k+1$ in $[k_{a,N}, k_{b,N}]$ we have a constant $C_{0}$ no depending from $k$ and $N$ such that  
$\vert \tilde \lambda^{(k+1)}_{N} - \tilde \lambda^{(k)}_{N}\vert \ge \frac {C_{0}}{N+2}$. 
\end{lemme}
\begin{preuve}{}
Since $f \in \mathcal C^4 ([0, 2 \pi])$, and $\rho \in \mathcal C^2 ([\theta_{1},\theta_{2}])$ we can write, using the main value theorem, 
$\tilde \lambda^{(k+1)}_{N} - \tilde \lambda^{(k)}_{N} = \frac{1}{N+2}  f'(c_{1})+ \frac{1}{(N+2)^2} \left( f'\rho\right)' (c_{2})
+ \frac{1}{(N+2)^3} \left(f'\rho'\rho\right)'(c_{3}) +\frac{1}{2}  \frac{1}{(N+2)^3} \left( f^{\prime\prime} \rho^2\right)' (c_{4})$
with $c_{j}\in ] \frac{k}{N+2}, \frac{k+1}{N+2}[$ for $j=1,2,3,4$. Then we obtain two constants $ \tau_{0}, M_{0}$, 
depending only from $f$ and $\rho$, 
such that $\vert \tilde \lambda^{(k+1)}_{N} - \tilde \lambda^{(k)}_{N}\vert \ge \frac{\tau_{0}}{N+2} -  \frac{M_{0}}{(N+2)^2}$.
Hence if $C_{0} = \frac{\tau_{0}}{N+2}$ we obtain  the result for a sufficiently  large $N$. 
\end{preuve}
\begin{lemme} \label{L22}
For $k,k+1$ in $[k_{a,N}, k_{b,N}]$ we consider the two eigenvalues $\lambda_{k} = \tilde \lambda^{(k)}_{N} +\frac{R_{1}}{(N+2)^3}$ and 
$\lambda_{k+1} = \tilde \lambda^{(k+1)}_{N} +\frac{R_{2}}{(N+2)^3}$. Then we have a constant $C_{1}$ no depending from $k$ such that for a sufficiently large $N$ we have $ \vert \lambda_{k+1} - \lambda_{k} \vert \ge \frac{C_{1} }{N+2}$.
\end{lemme}
\begin{preuve}{}
The property \ref{beautiful} implies than the rest in the equation (\ref{DOUZE}) is bounded by $\frac{T}{(N+2)^3}$ where $T$ is a real no depending from $\lambda$ . Hence 
with Lemma (\ref{L2}) we can write 
$ \vert \lambda_{k+1} - \lambda_{k} \vert \ge \frac{C_{1}}{N+2} - \frac{T}{(N+2)^3}$ and  with, for instance, 
$C_{1} = \frac{C_{0}}{2}$ we obtain the lemma for a sufficiently large $N$.
\end{preuve}
\begin{lemme}\label{L3}
For a fixed $k$ the equation (\ref{NEUF}) has one and only one solution in $[\theta_{1},\theta_{2}] $.
\end{lemme} 
\begin{preuve}{}
Assume $\tilde \lambda_{N}$ and $\tilde \lambda'_{N}$ two solutions 
of (\ref{NEUF}) for a same integer $k$. By (\ref{DOUZE}) and Property \ref{beautiful}
we have $\vert \tilde \lambda_{N}-\tilde \lambda_{N}'\vert \le \frac{T'}{(N+2)^3})$ with $T'$ no depending from $\lambda$.
By \cite{BARBA} we know that we have an eigenvalue
$\lambda_{N+1}$ of the matrix $T_{N+1}(f)$ with the bound $\tilde \lambda_{N} < \lambda_{N+1}
< \tilde \lambda'_{N}$ that implies 
$ \vert\tilde  \lambda_{N}- \lambda_{N+1}\vert =\le \frac{T'}{(N+2)^3})$. By Lemma \ref{L2} we have
$\vert \tilde \lambda_{N}- \lambda_{N+1}\ \vert \ge (\frac{C_{1}}{N+2})$,
that is a contradiction with the previous estimation. 
\end{preuve}
\section{Proof of Theorem \ref{THEO3}}
First we can observe that we can define the function $H$ on $[0, \pi] \times [0, \pi] $ with 
$H(0,0) =. f^{\prime \prime} (0)$ and  $H(\pi,\pi)= f^{\prime \prime} (\pi,\pi)$.
Hence for all $\lambda\in I_{0,2\pi}$ we have 
 $$ f(\theta) - \lambda=f_{1}(1-\cos \theta ) -\lambda=\left((1-\cos \theta)-(1-\cos \theta_{\lambda})\right) 
 H_{\lambda}(\theta)$$
 where $H_{\lambda} : \theta \mapsto H(\theta, \theta_{\lambda})$ is a regular function on $[-\pi, \pi]$.
  for all $\lambda \in [0, \pi].$
  With the same notations as in the proof of Theorem \ref{THEO1ter}  we can still write, for $\lambda \notin \operatorname{Spec}\left(T_{N}(f)\right)$
   $$T_{1,N,\lambda} = \frac{\det \left(T_{N-1}(f) -\lambda I_{N-1}\right)}
 {\det \left(T_{N}(f) -\lambda I_{N}\right)}.$$
As previously we have  also  $$\lambda \in \operatorname{Spec} \left(T_{N}(f)\right) \iff \frac{\det \left(T_{N}(f) -\lambda I_{N-1}\right)}
 {\det \left(T_{N-1}(f) -\lambda I_{N}\right)}=0.$$
 and, always with the equation (\ref{lastone}), we can write 
  \begin{equation} \label{INV11TER}
\frac{\det \left(T_{N}(f) -\lambda I_{N-1}\right)}
 {\det \left(T_{N-1}(f) -\lambda I_{N}\right)} = \frac{1-\bar \chi_{\lambda} ^{2(N+1)} \tau_{N}(\chi_{\lambda})}
 {\left(1-\bar \chi_{\lambda} ^{2(N+2)} \tau_{N}(\chi_{\lambda}) \right ) B_{2,N,\lambda} -B_{1,N,\lambda}},
 \end{equation}
 with 
 $$ \tau_{N}(\theta_{\lambda}) = 
 \frac{\bar P_{N+1,\lambda} (\chi_{\lambda}) P_{N+1,\lambda} (\chi_{\lambda})}
 {\bar P_{N+1,\lambda} (\overline{\chi_{\lambda}}) P_{N+1,\lambda} (\overline{\chi_{\lambda}})},$$
 and 
 $ B_{1,N,\lambda},B_{2,N,\lambda}$ as previously.  Hence we can write
\begin{equation}\label{SIX10}
 \lambda \in \left(\lambda \in \operatorname{Spec} \left(T_{N}(f)\right)\right) \cap I_{0,\pi}\iff
  \chi_{\lambda}^{2(N+2)}= \tau_{N}(\theta_{\lambda}),\lambda \in I_{0,\pi}.
  \end{equation}
  Since the function $H_{\lambda}$ is even, the constant 
$\tau_{N} (\theta_{\lambda})$ can be rewritten as
$$\tau_{N} (\theta_{\lambda}) = \left(\frac{P_{N+1,\lambda} ( \chi_{\lambda})}
{P_{N+1,\lambda} (\bar \chi_{\lambda})}\right)^2.
$$
On the other hand the function 
  $ \theta \mapsto \frac{P_{N+1,f(\theta)} (\bar \chi_{\lambda})}
{P_{N+1,f(\theta)} ( \chi_{\lambda})}$ is continuous from $[0,2 \pi]$ to $\{z \vert \vert z\vert =1\}$ hence we have a function $\rho_{N}$ defined and continuous on $[0,\pi]$        
  such that 
  $\tau_{N}(\theta_{\lambda}) = e^{2 i \rho_{N}(\theta_{\lambda})}$.
  Then equation (\ref{SIX10}) can be written 
    \begin{equation}\label{SEPT10}
 \lambda \in \left(\lambda \in \operatorname{Spec} \left(T_{N}(f)\right)\right) \cap I_{0,\pi}  \iff
\theta_{\lambda} = \frac{\rho_{N} (\theta_{\lambda}) +  k \pi}{ (N+ 2)} , k \in[0, 2N+3].
  \end{equation}
  Hence for $ k$ in $\{0, \cdots, 2N+3\}$  we have to find the solution in $[0, \pi]$ of the equation 
   \begin{equation}\label{HUIT10}
(N+2)\theta - \rho_{N} (\theta) =  k \pi
\end{equation}
But in the particular case where $\theta_{1} = 0$ and $\theta_{2}=\pi$ it is easy to verify that 
the function $\rho_{N}$ is in fact an odd $2\pi$-periodic function with $\rho_{N}(0)=0$.\\
Now we denote by $F_{N} $ the function $\theta\mapsto (N+2)\theta- \rho_{N}(\theta)$. For $1\le k \le N+1$ we have 
$$ F_{N}(0) = 0 < \pi k,  F_{N}(\pi) = (N+2) \pi > \pi k. $$
Hence the equation (\ref{HUIT10}) has at less one solution in $[0,\pi]$ for all $k \in \{ 0, \cdots N+1\}$.
In the other hand it is obvious that the solution of the equations
$\theta = \frac{\rho_{N} (\theta) +  k \pi}{ (N+ 2)}$ and $\theta = \frac{\rho_{N} (\theta) +  k' \pi}{ (N+ 2)}$
are different for $k\neq k'$. 
Since $f$ is strictly increasing on $[0, \pi]$, we have found $N+1$ eigenvalues of $T_{N}(f)$ in $[f(0), f(\pi)]$, and we will not obtain other 
eigenvalues outside the set $ \{ 0, \cdots N+1\}$.
The rest of the proof is the same as the proof of Theorem \ref{THEO1ter}.

  \section{Appendix}
 For the proof of Theorem \ref{THEO1ter} we have to know 
 $T_{N}(f)^{-1} _{1,1}$. 
 First we use Theorem \ref{THEO2} to obtain  $T_{N}(f_{r})^{-1} _{1,1}$ 
 with 
 $ f_{r}= \chi_{\lambda} (1-r \bar \chi_{\lambda}\chi)  (1-r \bar \chi_{\lambda}\bar \chi) \frac{1}{\vert P_{N+1,\lambda}\vert ^2}$, and now 
 $g_{1}=  \chi_{\lambda} (1-r \bar \chi_{\lambda}\chi) \frac{1}{P_{N+1,\lambda}}$, 
 $ g_{2}= (1-r \bar \chi_{\lambda}\bar \chi) \frac{1}{\overline{P_{N+1,\lambda}}}$\\
 We have to observe that $T_{N}(f_{r})^{-1} _{1,1}$ that is also $\langle T_{N}(f_{r})^{-1} (1) \vert 1\rangle$.
 Write $\langle T_{N}(f_{r})^{-1} (1) \vert 1\rangle = x_{0}-y_{0}$.
 Theorem \ref{THEO2} provides 
 $$ x_{0} = \langle \pi_{+} \left( \frac{1}{g_{2}}\right) \vert \frac{1}{\bar g_{1}}\rangle
 =  \frac{ \chi_{\lambda} }{ \vert P_{N+1,\lambda}(0)\vert ^{2}}.$$
 To obtain $y_{0}$ we need the terms 
 $\pi_{+} \left( \tilde \Phi_{N} \pi_{+}\left( \frac{1}{g_{2}}\right) \right)$ and 
 $\pi_{+} \left (\bar \Phi_{N} \pi_{+}\left( \frac{1}{\bar g_{1}}\right) \right)$.
 We have , if $\omega = r \bar \chi_{\lambda}$,
 $$ \pi_{+} \left (\tilde \Phi_{N} \pi_{+}\left( \frac{1}{g_{2}}\right) \right) =  
\overline {P_{N+1,\lambda}(0)} \pi_{+}\left( \frac{g_{2}}{g_{1}} \chi^{-N-1}\right)= 
C_{1} \frac{1}{1-\omega \chi} $$
with 
$$ C_{1}=\overline {P_{N+1,\lambda}(0)}\bar \chi_{\lambda}  \left( \frac{P_{N+1,\lambda}(\frac{1}{\omega} )}{\bar P_{N+1}(\omega)}\right) \omega^{N+1} (1-\omega^2).$$
Likewise we can write 
$$ \pi_{+} \left ( \bar \Phi_{N} \pi_{+}\left( \frac{1}{\bar g_{1}}\right) \right) 
= C'_{1} \frac{1}{1- \bar \omega \chi}, $$
with 
$$C'_{1}= \overline {P_{N+1,\lambda}(0)}\bar \chi_{\lambda} 
\left( \frac{P_{N+1,\lambda}(\frac{1}{\bar \omega} )}{\bar P_{N+1,\lambda}(\bar \omega)}\right) {\bar \omega}^{N+1} (1-{\bar \omega}^2).$$
  Hence 
  $$y_{0} =  C_{1}\overline{C'_{1}} \Bigl \langle (I-H_{\Phi_{N}^{*}} H_{\Phi_{N}} )^{-1}
  \frac{1}{1-\omega\chi} \big \vert  \frac{1}{1-\bar \omega\chi} \Bigr \rangle.$$ 
  We have now to use the following lemma 
  \begin{lemme}
  $\frac{1}{1-\omega \chi}$ is an eigenvector of $ H^\star _{\Phi_{N}} H _{\Phi_{N}} $ for the eigenvalue 
  $\tau_{N,r} (\omega)\omega^{2(N+2)} $ with \\
  $\tau_{N,r} (\omega)  =\frac{\overline{P_{N+1,\lambda}}(\frac{1}{\omega})
  P_{N+1,\lambda}(\frac{1}{\omega})}{\overline{P_{N+1,\lambda}}(\omega)P_{N+1,\lambda}(\omega)}$, with 
  $\Bigl\vert \omega^{2(N+2)}  \tau_{N,r} (\omega)\Bigr \vert < 1$ for $r\to1$ and $N$ sufficiently large.
  \end{lemme}
 It is Lemma 1 of \cite{RS020}. We obtain 
 $$ y_{0} =  C_{1}\overline{C'_{1}} \frac{1}{1-\omega^{2N+2} \tau_{N,r}(\omega)}  \frac{1}{1-\omega^{2}}.$$
 If now we consider the function $f_{1}$ defined by the  product $f_{1}= \tilde g_{1} \tilde g_{2}$ with 
 $\tilde g_{1} = \chi_{0}(1-\bar\chi_{0} \chi) \frac{1}{ P_{N+1}}$
 and $\tilde g_{2} = (1-\bar\chi_{0} \bar\chi) \frac{1}{\bar P_{N+1}}$,
 then for a fixed $N$, $\lim_{r\to 1} \left( T_{N}f_{r}\right)^{-1}_{1,1} = \left( T_{N}f\right)^{-1}_{1,1}$. Indeed 
 $$  \left( T_{N}f_{r}\right)^{-1}\left( T_{N}f\right) = \left( T_{N}f_{r}\right)^{-1}
 \left( T_{N}f_{r}\right) + \left( T_{N}f_{r}\right)^{-1} \left( T_{N}(f_{1}-f_{r})\right).$$
And $\lim_{r\to 1}  \left( T_{N}(f-f_{r})\right) =0$ that implies 
$ \lim_{r\to 1} \left( T_{N}f_{r}\right)^{-1}\left( T_{N}f\right) =I_{N}.$ 
Hence we can conclude that 
\begin{equation} \label{lastone}\left((T_{N}(f))^{-1}\right) _{1,1} = \frac{\left( 1-{\bar\chi_{\lambda}}^{2(N+2)}
\tau_{N}(\chi_{\lambda})\right ) B_{2,N,\lambda}- B_{1,N,\lambda}}
{1-\bar{\chi_{\lambda}}^{2(N+2)}\tau_{N} (\chi_{\lambda})},
\end{equation}
with 
$ B_{1,N,\lambda}=
C_{1}\overline{C'_{1}} (1-{\bar \chi_{\lambda}}^2)^{-1}$,
$B_{2,N,\lambda} =   \frac{ \chi_{\lambda}}{ \vert P_{N+1,\lambda}(0)\vert ^{2}}$,
and $\tau_{N}(\chi_{\lambda}) = \frac{\overline{P_{N+1,\lambda}}(\chi_{\lambda})
  P_{N+1,\lambda}( \chi_{\lambda})}{\overline{P_{N+1,\lambda}}(\bar\chi_{\lambda})P_{N+1,\lambda}(\bar\chi_{\lambda})}$

     \bibliography{Toeplitzdeux}

\end{document}